\documentclass{llncs}

\usepackage{times,epsfig,amssymb,latexsym,lscape}

\newcommand{\lb}{\hbox{{\it lb}}}
\newcommand{\rb}{\hbox{{\it rb}}}

\newcommand{\allEq}{\hbox{{\it allEq }}}

\title{Using propagation for solving complex arithmetic constraints}
\author{M.H. van Emden \and B. Moa}

\institute{Department of Computer Science, \\University of Victoria,
Victoria, Canada \\
\email{\{vanemden, bmoa\}@cs.uvic.ca},\\ WWW home page:
\texttt{http://www.cs.uvic.ca/\homedir vanemden/ }\\
April 4, 2003; as submitted to CP'03
}

\begin{document}
\maketitle

\begin{abstract}
Solving a system of nonlinear inequalities is an important problem for
which conventional numerical analysis has no satisfactory method.  With
a box-consistency algorithm one can compute a cover for the solution set
to arbitrarily close approximation.  Because of difficulties in the use
of propagation for complex arithmetic expressions, box consistency is
computed with interval arithmetic.  In this paper we present theorems
that support a simple modification of propagation that allows complex
arithmetic expressions to be handled efficiently. The version of box
consistency that is obtained in this way is stronger than when interval
arithmetic is used.

\end{abstract}

\section{Introduction}
One of the most important applications of constraint programming is
to solve a system of numeric inequalities:

\begin{equation}
\label{nonLinSys}
\begin{array}{ccccccccc}
g_1(x_1&,& x_2 &,& \ldots &,& x_m)  & \leq & 0     \\
g_2(x_1 &,& x_2 &,& \ldots &,& x_m) & \leq & 0     \\
\multicolumn{9}{c}{\dotfill}            \\
g_k(x_1 &,& x_2 &,& \ldots &,& x_m) & \leq & 0     \\
\end{array}
\end{equation}

Such  systems often appear as conditions in optimization problems.
Inequalities are regarded as intractable in conventional numerical
analysis.  The Kuhn-Tucker conditions allow these to be converted
to equalities.  The continuation method \cite{mrgn87a} is a fairly,
but not totally, dependable method for solving such
equalities. Moreover, it is restricted to polynomials.

An important contribution of constraint programming is the
box-consistency method \cite{bmcv94,vhlmyd97}, which improves on the
continuation method in several ways.  It applies not only to polynomials
$g_1, \ldots, g_k$ but to any functions that can be defined by an
expression that can be evaluated in interval arithmetic. It computes a
cover for the set of solutions and can make it fit arbitrarily closely.
In this way, all solutions are found and are approximated as closely as
required. The performance of the box-consistency method compares
favorably with that of the continuation method \cite{vhlmyd97}.

\section{Preliminaries}
In this section we provide background by reviewing some basic concepts.

\subsection{Constraint satisfaction problems}

In a {\em constraint satisfaction problem (CSP)} one has a set of
{\em constraints}, each of which is an instance of a formula.  Each of the
variables in the formula is associated with a {\em domain}, which is
the set of values that are possible for the variable concerned. A {\em
solution} is a choice of a domain element for each variable that makes
all constraints true.

With each type of constraint, there is an associated {\em domain
reduction operator}\/; DRO for short.  This operator may remove from the domains
of each of the variables in the constraint certain values that do not
satisfy the constraint, given that the other variables of the
constraint are restricted to their associated domains.
Any DRO is contracting, monotonic, and idempotent.

When the DROs of the constraints are applied in a ``fair'' order,
the domains converge
to a limit or one of the domains becomes empty.
A sequence of DROs activations is fair if every one of them occurs
an infinite number of times \cite{vnmd97,aptEssence}.
The resulting cartesian product of the domains becomes the
greatest common fixpoint of the DROs \cite{vnmd97,aptEssence}.
If one of the domains becomes empty,
it follows that no solutions exist within the initial
domains.  This, or any variant that leads to the same result, is called
a {\em constraint propagation algorithm}.

In practice, we are restricted to the domains that are representable in
a computer. As there are only a finite number of these,
constraint propagation terminates.

%A {\em constraint system} has the following components.
%\begin{enumerate}
%\item
%A set $\T = \{T_1, \ldots, T_m\}$ of sets called {\em types}.
%
%\item
%A set $\X = \{x_1, \ldots, x_n\}$ of {\em variables}.
%
%\item
%A {\em typing function} $\tau$: $\X \rightarrow \T$ such that for all
%$i\in\{1, \ldots, n\}$, $x_i$ is of type $\tau(x_i)\in \T$.
%
%\item
%A set of {\em constraints} $\C = \{C_1, \ldots, C_m\}$.  Each
%constraint is an atomic formula of first order predicate logic.
%
%\item
%A finite sequence of {\em domains} $\D = \langle D_1, \ldots, D_n
%\rangle$, where $D_i \subseteq \tau(x_i)$ for $i \in \{1, \ldots, n\}$.
%
%\end{enumerate}

%With each constraint $C$, we associate a so-called {\em Domain
%Reduction Operator} (DRO for short) that reduces the domain of each
%variable involved in $C$ by eliminating values that do not satisfy
%$C$.  That is, DRO contracts domains by removing inconsistent values.

\subsection{Constraint propagation}
A {\em Generic Propagation Algorithm} (GPA in the sequel) is a fair
sequence of DROs. A GPA maintains a pool of DROs, called {\em active set},
that still need to be applied. No order is specified for applying
these operators. Even though many variants of GPA exist
(see Apt~\cite{aptEssence} and Bartak~\cite{rthpra01}), they are all similar to
the pseudo-code given in Figure~\ref{lGPA}. Notice that the active
set $A$ is initialized to contain all constraints.

\begin{figure}
\begin{tabbing}
put all constraints into the active set $A$\\
while \=( $A \neq \emptyset$) $\{$\\
         \>choose a constraint $C$ from the active set $A$\\
         \>apply the DRO associated with $C$\\
         \> if one of the domains has become empty then stop\\
         \>add to $A$ all constraints involving variables whose
           domains have changed, if any\\
         \>remove $C$ from $A$\\
$\}$
\end{tabbing}
\caption{
\label{lGPA}
A pseudo-code for GPA.
}
\end{figure}

\subsection{Interval constraint satisfaction problems}
The constraint programming paradigm is very general. It applies to
domains as different as booleans, integers, finite symbolic domains,
and reals.  In this paper we consider {\em interval CSPs}, which are
CSPs where there is only one type and it is equal to the set
$\mathcal{R}$ of real numbers.  In such CSPs domains are restricted to
intervals, as reviewed below.

\subsection{Intervals}
A {\em floating-point number} is any element of $F \cup \{-\infty,
+\infty\}$,
where $F$ is a finite set of reals. A {\em floating-point interval} is
a closed connected set of reals, where the bounds, in so far as they exist,
are floating-point numbers. When we write ``interval'' without
qualification,
we mean a floating-point interval.
A {\em canonical interval} is a non-empty interval that does not properly
contain an interval. For every finite non-empty interval $X$, $lb(X)$ and
$rb(X)$ denote the left and right bound of $X$ respectively.
For an unbounded $X$, the left or right bound is defined as $-\infty$
or $+\infty$, or both.
Thus,
$ X = [lb(X), rb(X)]$ is a convenient notation for all non-empty intervals,
bounded or not.

If a real $x$ is not a floating-point number, then there is a unique
canonical interval containing it. Otherwise, there are three. Either way,
there is a unique least canonical interval containing $x$, and it is
denoted $\overline{x}$.

A {\em box} is a cartesian product of intervals.

\subsection{Box consistency}

In \cite{vhlmyd97}, box consistency is computed by a relaxation
algorithm implemented in interval arithmetic.  The algorithm takes as
input certain intervals $X_1, \ldots, X_m$ for the variables $x_1,
\ldots, x_m$.  It uses each of the functions $g_1, \ldots, g_k$ in the
way that is explained by a generic instance that we temporarily call
$g$.  We assume that the function $g$ is defined by an expression $E$
containing no variables other than $x_1, \ldots, x_n$.  We call the
algorithm in \cite{vhlmyd97} a {\em relaxation algorithm} because it
improves the intervals for the variables one at a time while keeping
the intervals for all the other variables fixed.
This is similar to the relaxation algorithms of numerical analysis.

For simplicity of
notation, let us assume that the interval for $x_1$ is to be improved
on the basis of the fixed intervals $X_2, \ldots, X_m$ for the
variables $x_2, \ldots, x_m$.  This is done by means of a function
$g_{X_2, \ldots, X_m}(x)$ that is defined by evaluating in interval
arithmetic the expression $E$ with $\overline{x}$ substituted for $x_1$
and $X_2, \ldots, X_m$ substituted for $x_2, \ldots, x_m$,
respectively.  Thus, $g_{X_2, \ldots, X_m}$ maps a real to an
interval.

To improve the interval $X_1 = [\lb(X_1),\rb(X_1)]$
for $x_1$, suppose that for some $a < \rb(X_1)$ we have
that
\begin{equation}
\label{trial}
\lb(g_{X_2, \ldots, X_m}([a,\rb(X_1)])) > 0.
\end{equation}
In that case the interval for $x_1$ can be improved from
$X_1$ to $[\lb(X_1),a]$.

A bisection algorithm is used to find the least floating-point $a$ for
which (\ref{trial}) holds, for fixed intervals $X_2, \ldots, X_m$.  A
similar bisection is used to improve the lower bound of $X_1$ using $g$
and the fixed intervals $X_2, \ldots, X_m$.  This exhausts what can be
done with $g$ and $X_2, \ldots, X_m$.  In general, repeating this
process with the other arguments and with the other functions among
$\{g_1,\ldots,g_k\}$ causes reductions of $X_2, \ldots, X_m$ and
further reductions of $X_1$.

If one of the intervals becomes empty, it is shown that no solution
exists within the original intervals $X_1, \ldots, X_m$. Otherwise, the
box-consistency algorithm terminates with no interval reduction possible.
Let us call the resulting state of the domains
{\em functional box consistency}, to
distinguish it from the type of box consistency described below.

As the criterion for $a$ being an improved upper bound for $X_1$ is
(\ref{trial})
with the left-hand side
evaluated in interval arithmetic, this box consistency algorithm can be
improved by means of interval constraints, as was pointed out in
\cite{vnmdn01b}.  Here it was proposed that instead of such interval
arithmetic evaluation, one applies propagation to the interval
CSP containing as constraints
\begin{eqnarray}
g(x_1,\ldots,x_m) & \leq & 0       \nonumber         \\
x_1 & > & a                                \\
x_2 \in X_2, \ldots, x_m \in X_m & & \nonumber
\end{eqnarray}
The result is called {\em relational box consistency}. In \cite{vnmdn01b}
it is shown that the resulting intervals are contained in those obtained by
functional box consistency.

To apply propagation, one needs to decompose the arbitrarily complex
expression for $g$ into multiple primitive arithmetic constraints
($x+y=z$, $x \times y=z$, as well as those involving trigonometric or logarithmic
constraints), as
explained in section~\ref{arch}, so that the corresponding efficient
DROs can be applied. In this way the structure of $E$
is lost. As a result, GPA
will activate DROs that have no effect,
even though, on eventual termination of GPA, the result is the correct one:
the unique and consistent state, or failure.

This problem was addressed in \cite{bggp99,vnmdn01b}. The remedy
described there is to create a tree data structure for $E$ and perform
propagation based on the tree structure. Such a {\em structured
propagation} algorithm avoids superfluous activations of DROs by
following the tree from the bottom to the top and then from the top to
the bottom and to repeat these two traversals as a cycle.

In this paper we show that by a simple modification of the GPA one can
get the effect of an optimized version of a structured propagation
algorithm, yet without maintaining a tree data structure.  The downward
part of this algorithm is enhanced by the absence of multiple occurrences
of variables.  In interval arithmetic this would be an unacceptable
limitation. However constraint programming allows us to translate
multiple occurrences to equality constraints.  This may be a welcome
improvement compared to the usual way of representing a system such as
in Equation~\ref{nonLinSys}.

In section~\ref{arch}, we describe our translation of a system such as
in Equation~\ref{nonLinSys} to a CSP consisting of primitive
constraints.

\section{Generating a CSP from a system of nonlinear numerical
inequalities}
\label{arch}

The functions $g_1, \ldots, g_k$ in Equation~\ref{nonLinSys}
are defined by expressions that can be evaluated in interval arithmetic
and hence can be translated to CSPs containing only primitive
constraints.  Any of these expressions may have multiple occurrences of
some of the variables.  As there are certain advantages in avoiding
multiple occurrences of variables in the same expression, we rewrite
without loss of generality the system in Equation~\ref{nonLinSys} to
the system shown in Figure~\ref{singleSys}.

\begin{figure}
\begin{displaymath}
\begin{array}{ccccccccc}
g_1(x_1&,& x_2 &,& \ldots &,& x_m) & \leq &0     \\
g_2(x_1 &,& x_2 &,& \ldots &,& x_m) & \leq & 0     \\
\multicolumn{9}{c}{\dotfill}            \\
g_k(x_1 &,& x_2 &,& \ldots &,& x_m) & \leq & 0     \\

\allEq(v_{1,1}&,& v_{1,2}&,& \ldots&,& v_{1,{n_1}}) &&  \\
\multicolumn{7}{c}{\dotfill}            \\
\allEq(v_{p,1}&,& v_{p,2}&,& \ldots&,& v_{p,{n_p}}) &&  \\
\end{array}
\end{displaymath}
\caption{
\label{singleSys}
A system of non-linear inequalities without multiple occurrences.
The variables $\{x_1,\ldots,x_m\}$ are partitioned into equivalence
classes
$V_1,\ldots,V_p$ where
$V_j$ is a subset $\{v_{j,1}, \ldots, v_{j,{n_j}}\}$
of $\{x_1,\ldots,x_m\}$,
for $j \in \{1,\ldots, p\}$.
}
\end{figure}

In Figure~\ref{singleSys}, the expressions for the functions $g_1,
\ldots, g_k$ have no multiple occurrences of variables.  As a result,
they have up to $m$ rather than $n$ variables, where $m \geq n$.  This
special form is obtained by associating with each of the variables
$x_j$ in Equation~\ref{nonLinSys} an equivalence class of the variables
in Figure~\ref{singleSys}.  In each expression each occurrence of a
variable is replaced by a different element of the corresponding
equivalence class. This can be done by making each equivalence class as
large as the largest number of multiple occurrences.
The predicate $\allEq$ is true if and only if all its real-valued
arguments are equal.

An advantage of this translation is that evaluation in interval
arithmetic of each expression gives the best possible result, namely
the range of the function values. At the same time, the $\allEq$
constraint is easy to enforce by making all intervals of the variables
in the constraint equal to their common intersection.  This takes
information into account from all $k$ expressions.  If the system in
its original form~\ref{nonLinSys}, with multiple occurrences, would be
translated to a CSP, then only multiple occurrences in a single
expression would be exploited at one time.

The translation of a complex arithmetic expression to a CSP containing
primitive arithmetic constraints is an obvious variant of the procedure
that has been familiar since FORTRAN compilers parsed a complex arithmetic
expression and generated code from the parse.  The CSP variant of this
procedure is at least as old as BNR Prolog \cite{BNR88},
which was first implemented
in the late 1980s. For a formal description of the translation
we refer to
\cite{vnmdn01b}.  We give a brief informal description here.

For the purpose of the translation, one regards an
expression as a tree with operators as internal
nodes and constants or variables as external nodes. We associate with
each internal node a unique variable.  Each internal node now generates
a primitive constraint. For example, if an internal node has ``$/$'' as
an operator, $x$ as a variable, and $y$ and $z$ as variables associated
with left and right child nodes, respectively, then the ternary
constraint $x \times z=y$ is generated.

In this way, each expression generates a primitive constraint for each
internal node.  Let $v$ be the variable associated with the root. As
this represents the value of the entire expression, the primitive
constraint $v \leq 0$ is generated as well. Because of the absence of
multiple occurrences, the partial CSP generated is an acyclic CSP for
which local consistency implies global consistency \cite{hyvo94}.

Finally, as the $\allEq$ constraints in Figure~\ref{singleSys} do not
contain expressions, they need no translation to more primitive
constraints: they have the obvious, optimal, and efficient domain
reduction operator described earlier.

This completes our description of how to translate a system as in
Equation~\ref{nonLinSys} to a CSP. To solve it
by propagation, we must first
consider propagation in a CSP consisting only of primitive constraints
generated from the tree of a single expression together with the
inequality constraint involving the root variable. This we do in the
next section.

Before proceeding thus, we point out that
translating the system (\ref{singleSys}) to a CSP in the way just described
enhances the opportunities for parallelism in propagation beyond those
already present in the system (\ref{nonLinSys}).
To investigate these would take us beyond the limits of this paper.
However, it will be useful here to highlight the structure that gives rise
to these opportunities by means of a hardware metaphor.

In the first place it is important to note that the sets of
constraints arising from the same expression
form a cluster for the purposes of propagation. For example, with the
exception of the root, the internal variables only have unique
occurrences.
As a result, when the DRO of a constraint is activated, it usually causes
constraints generated by the same expression to be added to the active
set. However, the external variables may occur in all of these clusters.

For the purposes of a parallel algorithm it is useful to imagine the
clusters arising from each of the expression as hardware ``cards'', each
connected to a ``bus'', where the lines of the bus represent the external
variables in common to several expressions.
%In figure~\ref{archfig} we have
%sketched such a hardware architecture.
%
%\begin{figure}[!htbp]
%\begin{center}
%\epsfxsize=5in
%\leavevmode
%\epsfbox{arch.eps}
%\caption{
%\label{archfig}
%A hardware architecture for CSP (2).
%}
%\end{center}
%\end{figure}
%
%The cards have been given a triangular
%form to symbolize the fact that they have the structure of expression
%trees.
Whenever two external variables belong to the same equivalence
class, they are connected by a ``jumper'' in the hardware model.

The hardware architecture suggests a parallel process for each card
that asynchronously executes DROs of constraints only involving
internal variables. The processes synchronize when they access one of
the bus variables. The DROs of $\allEq$ constraints can also be
executed by a parallel process dedicated to each.

\section{Modifying propagation for evaluating an expression}

We first show that, regardless of efficiency,
GPA can be used to evaluate a single expression.
Suppose that we have an expression $E$ in variables $x_1,\ldots,x_n$.
Let $y$ be the variable at the root of the tree representing
the value of $E$.
Let $C$ be the CSP generated by $E$ as described before.

\begin{theorem}\label{basic}
Suppose the domains of $x_1,\ldots,x_n$ are the intervals
$X_1,\ldots,X_n$.  Suppose the domains of $y$ and the other internal
variables are $[-\infty,+\infty]$.  Applying the GPA to $C$ results in
the domain of $y$ being the same interval as the one obtained by
evaluating $E$ in interval arithmetic with $X_1,\ldots,X_n$ substituted
for $x_1,\ldots,x_n$, respectively.

\end{theorem}

{\em Proof.}
According to \cite{vnmd97,aptEssence}, every fair
sequence of DROs in GPA converges to the same limit
for the domains of the variables. There
is a finite sequence $s$ of DROs that mimics the evaluation of $E$ in
interval arithmetic. At the end of this, $y$ has the value computed by
interval arithmetic and GPA terminates.

\paragraph{}
Theorem~\ref{basic} is useful in showing that, in the absence of
information about the value of the expression, propagation does the
equivalent of interval arithmetic.  But the GPA does it in a wasteful
way. GPA does not specify the order of applying the DROs other than
that their sequence should be a fair one.  In a typical random fair
sequence,
many DROs will not have any effect.  This inefficient
behavior is the motivation for our modifications
to GPA presented here.

\begin{theorem}
\label{modEval}
Suppose we modify GPA so that the active set is initialized to contain
instead of all constraints only those containing at most one internal
variable. Suppose also that the active set is a queue in which the
constraints are initially ordered according to the level they occupy in
the expression tree, with those that are further away from the root
placed nearer to the front of the queue.  Then GPA terminates after
activating the DRO of every constraint at most once.  On termination,
$y$ has as domain the value that the expression has in interval
arithmetic.

\end{theorem}

Thus we see that GPA has exactly the right behavior for the evaluation
of an expression if only we initialize the active set with the right
selection of constraints. We call this {\em propagation with selective
initialization (PSI)}. In the sequel, the constraints that have at most
one internal variable are referred to as {\em peripheral constraints}.

\section{Selective initialization for obtaining box consistency}

Propagation, whether modified or not, obtains results that are at least
as strong, and typically stronger than, box consistency as described in
\cite{vhlmyd97}.  This can already be demonstrated when considering a
single expression $E$ in variables $x_1, \ldots, x_m$ with intervals
$X_1, \ldots, X_m$ as domains.  A step towards box consistency is to
reduce the interval for each of the variables separately. To simplify
notation we do this for $x_1$, keeping the interval domains $X_2,
\ldots, X_m$ for $x_2, \ldots, x_m$ fixed.  Suppose that for some $a <
\lb(X_1)$ we have that $\lb(g_1([a,\rb(X_1)],X_2, \ldots, X_m)) > 0$.
In that case the interval for $x_1$ can be improved from $X_1$ to
$[\lb(X_1),a]$.

Suppose that instead of such an interval arithmetic evaluation, one
applies propagation to the interval constraint system containing as
constraints $x_1 > a$, $x_2 \in X_2, \ldots, x_m \in X_m$, as well as
all the ones obtained by translating the expression tree of $E$ to
primitive constraints.

\begin{theorem} \label{unModBox}
Suppose that $\lb(E([a,\rb(X_1)],X_2, \ldots, x_m \in X_m) >0$.
When GPA is applied to this CSP, failure results.

\end{theorem}

{\em Proof.} Consider any fair sequence that starts with a segment $s$
mimicking the interval arithmetic evaluation of $E$.  At the end of
$s$, the interval for $y$ has a positive lower bound, by the
assumption.  The fair sequence can be continued by applying the DRO for
$y \leq 0$.  This yields failure.

This proves the theorem.
If we use GPA in this
way for box consistency instead of interval arithmetic, we never obtain
a worse result and we typically obtain a better result.

However, unmodified GPA will obtain the better result in an inefficient
way.

Our next result is a modification of GPA that obtains the better result
in a more efficient way.

\paragraph{}

Suppose we have a CSP $S$ generated by an expression $E$. Let
$y$ be the variable at the root of the tree representing $E$. Suppose
we apply GPA to $S$. After the termination of GPA, we have the following
theorem.

\begin{theorem}
\label{modSolve}
Suppose the domain for $y$ is changed to a proper subset of it.
Suppose the unique
constraint containing $y$ is placed in the active set as only element.
Then
GPA terminates with the same result as when the active set would have
been initialized to contain all constraints.
\end{theorem}

\begin{proof}
Let $s$ be any fair sequence of DROs. If $s$ starts with the unique
constraint containing $y$, then the theorem follows. If $s$ starts with
a different constraint $c$, then applying the DRO associated with $c$
does not affect any domain (DROs are idempotent and the DRO of $c$ is
already in its fixpoint). Thus, removing $c$ from $s$ does not affect the
fixpoint of $s$.  We keep removing the first element of
$s$ until we reach the unique constraint containing $y$. The new $s'$
formed starts with the unique constraint containing $y$ and has the
same fixpoint as $s$. Thus, Theorem~\ref{modSolve} is proved.
\end{proof}

This theorem shows that a complex constraint of the form $E(x_1,
\ldots, x_n) \leq 0$ can be made
relationally box consistent by using a modification of GPA that is
more efficient without affecting the quality of the result.

The result obtained from
Theorem~\ref{modSolve} is usually better than
backward evaluation \cite{bggp99}. The following
example illustrates that.

Let us consider the following system.
\begin{eqnarray}
x_1 / x_2 & \leq & 0 \nonumber \\
 x_1 \in [-1,1] & & x_2 \in [-1,1] \nonumber
\end{eqnarray}
Using the decomposition described in section~\ref{arch},
we generate the following CSP.
\begin{eqnarray}
y & \leq & 0 \nonumber\\
x_1 / x_2 & = &y \nonumber\\
 x_1 \in [-1,1], & x_2 \in [-1,1], & y \in [-\infty,+\infty] \nonumber
\end{eqnarray}

Applying the GPA described in Theorem~\ref{modSolve} gives the
fixpoint $[-1,1] \times [-1,1] \times [-\infty,-1]$. The result
obtained from interval arithmetic backward evaluation is $[-1,1] \times
[-1,1] \times [-\infty, 0]$.

Pseudo-code for the PSI algorithm is given in Figure~\ref{PSIAlg}.

\begin{figure}
\begin{tabbing}
put only \textbf{peripheral} constraints into the active set $A$\\
while \=( $A \neq \emptyset$) $\{$\\
         \>choose a constraint $C$ from the active set $A$\\
         \>apply the DRO associated with $C$\\
         \> if one of the domains is empty then stop\\
         \>add to $A$ all constraints involving variables whose
	   domains have changed, if any\\
         \>remove $C$ from $A$\\
$\}$
\end{tabbing}
\caption{
\label{PSIAlg}
Pseudo-code for Propagation by Selective Initialization.
}
\end{figure}

In certain situations (fully described in \cite{mdmo03}),
Theorem~\ref{modSolveStrong}, stated below, can be used instead
of Theorem~\ref{modSolve}.

In addition to the assumptions of  Theorem~\ref{modSolve},
we suppose that the expression $E$ has no multiple occurrences
of any variable.

\begin{theorem}
\label{modSolveStrong}
Suppose the domain for $y$ is changed to a proper subset of it.
Suppose the unique constraint containing $y$ is placed in the active
set as only element.  Then GPA terminates after having activated the
DRO of each constraint at most once.

Moreover, the result is the same as when the active set would have been
initialized to contain all constraints.
\end{theorem}

This theorem is based on the fact that once the fixpoint is obtained,
reducing the domain of one variable may cause the other domains
to be reduced but not itself. As shown in
\cite{mdmo03}, this is true when DROs
have certain properties. For example, one could allow domains
to be the union of disjoint intervals, as in
the systems of Hyv{\"o}nen or Havens \cite{hyvo94,hssjv92}.
But when DROs are those described in this paper, reducing the domain of a
variable can affect its own domain as shown in the example above. Even
though $[-1,1] \times [-1,1] \times [-\infty,+\infty]$ is a fixpoint
of the CSP
\begin{eqnarray}
x_1 / x_2 & = &y \nonumber\\
 x_1 \in [-1,1], & x_2 \in [-1,1], & y \in [-\infty,+\infty] \nonumber
\end{eqnarray}
reducing the domain of $y$ to $[-\infty,0]$ leads to a further
reduced domain of $y$ that is equal to $[-\infty,-1]$.

\section{Conclusion}
In this paper, we presented a slight modification to propagation to
get the same or better results than structured
propagation. Such a modification is used in a box-consistency
algorithm to solve systems of non-linear inequalities but
it can also be used when solving other CSPs obtained by translating complex
expressions.
In fact, most CSPs of practical importance seem to be derived from complex
expressions.

\section*{Acknowledgments}
We acknowledge generous support by the
University of Victoria, the Natural Science and Engineering Research
Council NSERC, the Centrum voor Wiskunde en Informatica CWI, and the
Nederlandse Organisatie voor Wetenschappelijk Onderzoek NWO.
%\bibliography{sample}
%\bibliographystyle{plain}
%\bibliography{article,book,collection,proceedings,techreport,thesis,unpublished,misc}

\end{document}